\documentclass{tac}
\usepackage{xy}
\xyoption{all}

\def\mathcaldef#1{\expandafter\def\csname#1\endcsname{{\cal#1}}}

\def\qq{\quad\quad}
\def\qv{\qq ;\qq}
\def\iso{\,\cong\,}

\def\ot{\leftarrow}

\def\imp{\Rightarrow}

\def\la{\langle}
\def\ra{\rangle}
\def\adj{\dashv}
\def\op{^{\rm op}}
\def\ast{^*}
\def\st{^{-1}}
\def\ex{\exists}
\def\fa{\forall}
\def\comp{\pi_0}
\def\ov{\overline}
\def\ul{\underline}
\def\tm{\times}
\def\dm{\diamondsuit}
\def\sq{\sqcap\!\!\!\!\!\sqcup}
\def\sqq{\sqcap\!\!\!\!\sqcup}
\def\cd{\cdot}
\def\Si{\Sigma}
\def\wdg{\wedge}
\def\eps{\varepsilon}
\def\bs{\backslash}

\mathcaldef{C}
\mathcaldef{D}
\mathcaldef{U}
\mathcaldef{E}
\mathcaldef{M}
\mathcaldef{S}
\mathcaldef{P}
\mathcaldef{B}
\mathcaldef{V}

\mathbfdef{Set}
\mathbfdef{Cat}
\mathbfdef{Pos}
\mathbfdef{Grph}

\mathrmdef{id}

\def\hom{\ul{\rm hom}}
\def\meets{\ul{\rm meets}}
\def\nat{\ul{\rm nat}}
\def\ten{\ul{\rm ten}}
\def\lim{\ul{\rm lim}}
\def\colim{\ul{\rm colim}}
\def\eend{\ul{\rm end}}
\def\coend{\ul{\rm coend}}
\def\true{{\tt true}}
\def\false{{\tt false}}
\def\sup{{\rm sup}}
\def\inf{{\rm inf}}
\def\Sup{{\rm Sup}}
\def\Inf{{\rm Inf}}

\newtheorem{prop}{Proposition}
\newtheorem{corol}{Corollary}
\let\pf\proof
\let\epf\endproof
\def\eq{\begin{equation}}
\def\eeq{\end{equation}}

\author{Claudio Pisani}
\address{via Gioberti 86,\\ 10128 Torino, Italy.}
\title{A logic for categories}
\copyrightyear{2010}
\eaddress{pisclau@yahoo.it}

\begin{document}

\maketitle
\date{\today}

\begin{abstract}
We present a doctrinal approach to category theory, obtained by abstracting from
the indexed inclusion (via discrete fibrations and opfibrations) of left and of right actions
of $X\in\Cat$ in categories over $X$.
Namely, a ``weak temporal doctrine" consists essentially of two indexed functors with the same codomain
such that the induced functors have both left and right adjoints satisfying some exactness conditions,
in the spirit of categorical logic.

The derived logical rules include some adjunction-like laws
involving the truth-values-enriched hom and tensor functors, which 
condense several basic categorical properties and display a nice symmetry.
The symmetry becomes more apparent in the slightly stronger context of 
``temporal doctrines", which we initially treat and
which include as an instance the inclusion of lower and upper sets in the parts of a poset,
as well as the inclusion of left and right actions of a graph in the graphs over it.
\end{abstract}

\tableofcontents

\section{Introduction}
\label{intro}

Let $X$ be a set endowed with an equivalence relation $\sim$, and let $\V X$ be the poset
of closed parts, that is those subsets $V$ of $X$ such that $x\in V$ and $x\sim y$ implies $y\in V$.
A part $P\in\P X$ has both a ``closure" $\dm P$ and an ``interior" $\sq P$, that is the 
inclusion $i:\V X\to\P X$ has both a left and a right adjoint:
\[ \dm \adj i \adj \,\sqq : \P X \to \V X \]  
Thus the (co)reflection maps (inclusions) $\eps_P:i\sq P\to P$ and $\eta_P:P\to i\dm P$ induce 
bijections (between 0-elements or 1-elements sets): 
\[
\begin{array}{c}
\P X(iV,i\sq P) \\ \hline
\P X(iV,P) 
\end{array}
\qv
\begin{array}{c}
\P X(i\dm P,iV) \\ \hline
\P X(P,iV) 
\end{array}
\]
By taking $\sq(P\imp Q)$ as a $\V X$-enrichment of $\P X(P,Q)$, it turns out that the above 
adjunctions are also enriched in $\V X$ giving isomorphisms:
\[
\begin{array}{c}
\sq(iV \imp i\sq P) \\ \hline
\sq(iV \imp P) 
\end{array}
\qv
\begin{array}{c}
\sq(i\dm P \imp iV) \\ \hline
\sq(P \imp iV) 
\end{array}
\]
We also have the related laws:
\[
\begin{array}{c}
i\dm(iV \tm iW) \\ \hline
iV \tm iW 
\end{array}
\qv
\begin{array}{c}
i\sq(iV \imp iW) \\ \hline
iV \imp iW 
\end{array}
\qv
\begin{array}{c}
\dm(i\dm P \tm iV) \\ \hline
\dm(P \tm iV) 
\end{array}
\]
the first two of them saying roughly that closed parts are closed with respect to product (intersection)
and exponentiation (implication).
Given a groupoid $X$, the same laws hold for the inclusion of the actions of $X$ in the groupoids over $X$ 
(via ``covering groupoids").
The above situation will be placed in the proper general context in sections~\ref{adj} and~\ref{hyper},
where we develope some technical tools concerning enriched adjunctions and apply them to hyperdoctrines~\cite{law70}.

Now, let us drop the symmetry condition on $\sim$, that is suppose that $X$ is a poset; 
then we have the poset of lower-closed parts $\D X$ and that of upper-closed parts $\U X$.
Again, the inclusions $i:\D X\to \P X$ and $i':\U X\to \P X$ have both a left and a right adjoint:
\[
\dm \adj i \adj \,\sqq : \P X \to \D X \qv \dm' \adj i' \adj \,\sqq' : \P X \to \U X 
\]
While some of the above laws still hold ``on each side":
\eq     \label{int1c}
\begin{array}{c}
\sq(iV \imp i\sq P) \\ \hline
\sq(iV \imp P) 
\end{array}
\qv
\begin{array}{c}
\sq'(i'V \imp i'\sq' P) \\ \hline
\sq'(i'V \imp P) 
\end{array}
\eeq
\eq      \label{int1d}
\begin{array}{c}
i\dm(iV \tm iW) \\ \hline
iV \tm iW 
\end{array}
\qv
\begin{array}{c}
i'\dm'(i'V \tm i'W) \\ \hline
i'V \tm i'W 
\end{array}
\eeq
the other ones hold only in a mixed way:
\eq    \label{int1a}
\begin{array}{c} 
\sq'(i\dm P \imp iV) \\ \hline
\sq'(P \imp iV) 
\end{array}
\qv
\begin{array}{c}  
\sq(i'\dm' P \imp i'V) \\ \hline
\sq(P \imp i'V) 
\end{array}
\eeq
\eq   \label{int1}
\begin{array}{c}
i\sq(i'V \imp iW) \\ \hline
i'V \imp iW 
\end{array}
\qv
\begin{array}{c}
i'\sq'(iV \imp i'W) \\ \hline
iV \imp i'W 
\end{array}
\eeq
\eq     \label{int1b}
\begin{array}{c}
\dm(i\dm P \tm i'V) \\ \hline
\dm(P \tm i'V) 
\end{array}
\qv
\begin{array}{c}
\dm'(i'\dm' P \tm iV) \\ \hline
\dm'(P \tm iV) 
\end{array}
\eeq
The laws~(\ref{int1c}) through~(\ref{int1b}) hold also for the inclusion of the left and the right actions
of a category $X$ in categories over $X$ (via discrete fibrations and opfibrations): 
\[
i:\Set^{X\op} \to \Cat/X \qv i':\Set^X \to \Cat/X
\]
and, when they make sense, also for the inclusion of open and closed parts in the parts of a topological space
(or, more generally, of local homeomorphisms and proper maps to a space $X$ in spaces over $X$;
see~\cite{pis09}).
 
Abstracting from these situations, we may define a ``temporal algebra" as a cartesian closed category
with two reflective and coreflective full subcategories satisfying the above laws 
(in fact, it is enough to assume either~(\ref{int1a}) or~(\ref{int1}) or~(\ref{int1b})).
A ``temporal doctrine" is then essentially an indexed temporal algebra $\,\la i_X:\M X\to\P X\ot\M'X:i'_X\, ;\, X\in\C\ra\,$
such that the inclusions $i_1$ and $i'_1$ over the terminal object $1\in\C$ are isomorphic.
Temporal doctrines and their basic properties are presented in sections~\ref{temp} and~\ref{fp}.

In Section~\ref{tv} we show how the ``truth-values" $\M1\iso\M'1$ serve as values 
for an enriching of $\P X$, $\M X$ and $\M'X$ in which the adjunctions  
\[    
\Si_f \adj f\ast \adj \,\Pi_f : \P X \to \P Y 
\]
\[    
\dm_X \adj i_X \adj \,\sqq_X : \P X \to \M X  \qv \dm'_X \adj i'_X \adj \,\sqq'_X : \P X \to \M'X
\]
\eq   \label{int4}
\ex_f \adj f\cd \adj \,\fa_f : \M X \to \M Y  \qv  \ex'_f \adj f'\cd \adj \,\fa'_f : \M'X \to \M'Y 
\eeq
are also enriched (where, for $f:X\to Y$ in $\C$, $\ex_f M \iso \dm_Y\Si_f i_X M$).
For example, the temporal doctrine of posets is two-valued while that of reflexive graphs is $\Set$-valued,
by identifying sets with discrete graphs.
 
If $\C = \Cat$, the functors $\Pi_f$ are not always available, and the above mentioned enrichment is only
partially defined. 
This weaker situation will be axiomatized in Section~\ref{cat},
where we will see that~(\ref{int4}) still can be enriched giving:
\eq   \label{int5}
\begin{array}{c}
\nat_X(L,f\cd M) \\ \hline
\nat_Y(\ex_f L,M)
\end{array}
\qv  
\begin{array}{c}
\nat_X(f\cd M, L) \\ \hline
\nat_Y(M,\fa_f L)
\end{array}
\eeq
(and similarly for ``right actions" or ``right closed parts" in $\M'$) where 
\[
\nat_X(L,M) := \eend_X(i_X L \imp i_X M) := \fa_X\sq_X(i_X L \imp i_X M)
\]
In a somewhat dual way one also obtains:
\eq    \label{int6}
\begin{array}{c}
\ten_X(N,f\cd M) \\ \hline
\ten_Y(\ex'_f N,M)
\end{array}
\qv  
\begin{array}{c}
\ten_X(f'\cd N,L) \\ \hline
\ten_Y(N,\ex_f L)
\end{array}
\eeq
where one defines the tensor product by
\[
\ten_X(N,M) := \coend_X(i'_X N \tm i_X M) := \ex_X\dm_X(i'_X N \tm i_X M)
\]
In Section~\ref{lim} we show how the laws~(\ref{int5}) and~(\ref{int6}) allow one to derive 
in an effective and transparent way several basic facts of category theory, in particular concerning 
(co)limits, the Yoneda lemma, Kan extensions and final functors.
In section~\ref{compr} and~\ref{sup} other ``classical" properties are obtained exploiting
also a ``comprehension" axiom, relating $\P X$ and $\C/X$. 

This approach also offers a new perspective on duality: we do not assume the coexistence of a (generalized) category $X$
and its dual $X\op$ (which in fact is not so obvious as it may seem at a first sight).
Rather, we capture the interplay between left and right ``actions" or ``parts" of a category 
or ``space" $X$ by the above-sketched axioms concerning the inclusion of both of them in 
a category of more general ``labellings" or ``parts".

It is remarkable that while~(\ref{int1}) is equivalent to~(\ref{int1a}) and to~(\ref{int1b}),
they underlie seemingly unrelated items. 
On the one hand, for a truth value $V$ in $\M 1\iso\M'1$ the ``$V-$complement" of $M\in \M X$ 
\[
\neg(M,V) \,\, := \,\, i_X M \imp i'_X X'\cd V  
\]
``is valued" in $\M'X$, that is factors by~(\ref{int1}) through $i'_X$ (and conversely).
This generalizes the open-closed duality via complementation in topology
(and in particular the upper-lower-sets duality for a poset) which is given by $\neg(M,\false)$.

On the other hand, if we denote by $\{ x \} := \Si_x 1$ the ``part" in $\Cat/X$ corresponding
to the object $x:1\to X$, then $\dm_X\{ x \} = X/x$ corresponds to the presheaf represented by $x$
and for $N\in\M'X$ we can prove that $\ten_X(N,X/x) \iso x'\cd N$ using~(\ref{int1b}) as follows: 
\[
\ten_X(N,X/x) \iso \ex_X\dm_X(i'_X N \tm_X i_X\dm_X\{ x \}) \iso \ex_X\dm_X(i'_X N \tm_X \{ x \}) \iso 
\] 
\[
\iso \dm_1\Si_X(i'_X N \tm_X \Si_x 1) \iso \dm_1\Si_X\Si_x(x\ast i'_X N \tm 1) \iso \dm_1 i_1 (x'\cd N)  \iso x'\cd N
\] 
While in any temporal doctrine we can similarly derive $\hom_X(X/x,M) \iso x\cd M$ using~(\ref{int1a}), 
in $\Cat$ such a proof of the Yoneda lemma stumbles against the lack of $\Pi_x$ 
and the related non-exponentiability of $\{x\}$ (whenever $x$ is a non-trivial retract in $X$).
In this case, or in any (weak) temporal doctrine, one can use directly the first of~(\ref{int5}):
\[
\nat_X(X/x,M) \iso \nat_X(\ex_x 1 ,M) \iso \nat_1(1 , x\cd M) \iso x\cd M
\]
Similarly, using~(\ref{int6}) one gets again:    
\[
\ten_X(N,X/x) \iso \ten_X(N,\ex_x 1) \iso \ten_1(x'\cd N,1) \iso x'\cd N
\]

The present paper is a development of previous works on ``balanced category theory"
(see in particular~\cite{pis08} and~\cite{pis09}); 
the doctrinal approach adopted here emphasizes the logical aspects and suggests a wider range of applications.

\section{Enriching adjunctions}
\label{adj}

In this section, we make some remarks that will be used in the sequel. 
Along with ordinary adjunctions, Kan defined and studied what are now known as 
adjunctions with parameter and enriched adjunctions. 
In particular, we will use the following result from~\cite{kan}: 
\begin{lemma}   \label{par}
Given functors $F,F':\C\tm\P\to\V$ and $R,R':\P\op\tm\V\to\C$ such that there are 
adjunctions (with parameter) 
\[ 
F(X,P)\adj R(P,V) \qv F'(X,P)\adj R'(P,V)
\]
the natural transformations $F\to F'$ correspond bijectively to the ones $R'\to R$,
and this correspondence restricts to natural isomorphisms.
In particular $F\iso F'$ iff $R\iso R'$.
\epf
\end{lemma}
The next remark roughly says that a geometric morphism is naturally enriched in its codomain:
\begin{prop}   \label{prod}
Let $F\adj R:\C\to\V$ be an adjunction between cartesian closed categories, with $F$ left exact. 
Then $\C$ is enriched in $\V$ via 
\eq
\hom_\C^\V(X,Y) := R(X\imp_\C Y)
\eeq
and there are natural isomorphisms:
\[
\hom_\C^\V(FV,X) \iso \hom_\V(V,RX)
\]
(where $\hom_\V(V,W) := V\imp_\V W$ is the internal hom of $\V$)
that is the adjunction $F\adj R$ is itself enriched in $\V$.
Furthermore, the natural transformations given by the arrow mappings of $F$ and $R$ are also enriched: 
\[ 
\hom_\V(V,W) \to \hom_\C^\V(FV,FW)    \qv     \hom_\C^\V(X,Y) \to \hom_\V(RX,RY)
\]
\end{prop}
\pf
For the first part, we have
\[
\V(1,R(X\imp Y)) \iso \C(F1,X\imp Y) \iso \C(1,X\imp Y) \iso \C(X, Y) 
\]
(In fact, more generally, $R$ transfers any enriching in $\C$ to an enriching in $\V$.)
For the second part, since $F(V\tm W)\iso FV\tm FW$, we can apply Lemma~\ref{par}
to the adjunctions:
\[
F(V\tm W) \adj W\imp RX    \qv   FV \tm FW \adj R(FW\imp X)
\]
For the third part, the chain of natural transformations:
\[
\V(U,V\imp W)) \iso \V(U\tm V,W) \to \C(F(U\tm V),FW) \iso 
\]
\[
\iso \C(FU\tm FV,FW) \iso \C(FU,FV\imp FW) \iso \V(U,R(FV\imp FW))
\]
yields the desired natural transformation, which is easily seen to enrich the arrow mapping of $F$.
For $R$ we similarly have:
\[
\V(U,R(X\imp Y)) \iso \C(FU,X\imp Y) \iso \C(FU\tm X,Y) \to 
\]
\[
\to \C(F(U\tm RX),Y) \iso \V(U\tm RX,RY) \iso \V(U,RX\imp RY)
\]
where the non-isomorphic step is induced by the canonical 
\eq  \label{frob}
\la Fp,\eps Fq \ra : F(U\tm RX) \to FU\tm X
\eeq
Thus, $R$ is fully faithful, also as an enriched functor, iff~(\ref{frob}) is an iso,
that is $F\adj R$ satisfies the Frobenius law. 
Since here we have not used the fact that $F$ preserves all finite products, but only
the terminal object (in order to obtain an erichment of $R$) we get 
in particular a proof of Corollary 1.5.9 (i) in~\cite{ele}.   
\epf
If $F$ has a further left adjoint $L:\C\to\V$, then it is left exact and the above proposition applies.
We now show that in this case the adjunction $L\adj F$ is also enriched in $\V$ iff it satisfies the 
Frobenius reciprocity law:
\begin{prop}  \label{left}
Suppose that $\C$ and $\V$ are cartesian closed and that 
\[
L\adj F\adj R:\C\to\V
\]
Then the existence of the following natural isomorphisms are equivalent:
\begin{enumerate}
\item
\(
LX \tm_\V V \iso L(X\tm_\C F V)
\)
\item
\(
F(V\imp_\V W) \iso FV\imp_\C FW
\)
\item
\(
\hom_\C^\V(X,FV) \iso \hom_\V(LX,V)
\)
\end{enumerate}
\end{prop}
\pf
As before, we apply Lemma~\ref{par} to the adjunctions:
\[
LX \tm V  \adj F(V\imp W)  \qv   L(X\tm FV) \adj FV\imp RW
\]
getting the equivalence of 1) and 2) (which is well-known; see e.g.~\cite{law70}),
and to the adjunctions:
\[
LX \tm V  \adj LX\imp W  \qv   L(X\tm FV) \adj R(X\imp FW)
\]
getting the equivalence of 1) and 3).
\epf
Note that the same functor $\C\tm \V \to \V$ has two different right adjoints,
depending on the parameter chosen.

\begin{remark}   \label{ff}
It is well known that given adjunctions $L\adj F\adj R:\C\to\V$, with $F$ fully faithful,
if $\C$ is cartesian closed then so is also $\V$; in fact, products in $\V$ can be defined by
\[
1_\V := R1_\C \qv V\tm_\V W := R(FV\tm_\C FW) 
\]
or also by
\[
1_\V := L1_\C \qv V\tm_\V W := L(FV\tm_\C FW) 
\]
and exponentials by
\eq    \label{exp}
V\imp_\V W :=  R(FV\imp_\C FW)
\eeq
Note that, following Proposition~\ref{prod},~(\ref{exp}) indicates that
$F$ is fully faithful as an enriched functor, and we get
\eq       \label{ff1}
\hom_\C^\V(FV,X) \iso \hom_\C^\V(FV,FRX)
\eeq
Note also that, in this case, the equivalent conditions of Proposition~\ref{left} can be
rewritten as follows:
\eq   \label{ff2}
L(X\tm FV) \iso  L(FLX \tm FV) 
\eeq
\eq    \label{ff3}
FV\imp FW \iso FR(FV\imp FW) 
\eeq
\eq    \label{ff4}
\hom_\C^\V(X,FV) \iso \hom_\C^\V(FLX,FV)
\eeq
where the isomorphisms are induced by the unit of $L\adj F$ (the first and the third ones)
and by the counit of $F\adj R$ (the second one). 
\end{remark}

\section{The logic of hyperdoctrines}
\label{hyper}

We now show how some of the results of Section~\ref{adj} apply to hyperdoctrines~\cite{law70}, 
giving interesting consequences.
Recall that a hyperdoctrine is an indexed category $\,\la \P X\, ;\, X\in\C\ra\,$ such that $\C$
and all the categories $\P X$ are cartesian closed, and such that each substitution functor $f\ast:\P Y\to\P X$
has both a left and a right adjoint $\Si_f\adj f\ast\adj\Pi_f$ for any $f:X\to Y$ in $\C$.
The logical significance of hyperdoctrines, and in particular the
role of the adjoints to the substitution functors as existential and universal quantification, 
and that of $\P 1$ as ``sentences" or ``truth values",
are clearly illustrated in~\cite{law70} and in other papers by the same author.

Here we also assume that the adjunctions $\Si_f\adj f\ast$ satisfy the Frobenius law. 
On the other hand, we do not need to assume that $\C$ is cartesian closed but
only that it has a terminal object.
\begin{corollary}       \label{hyp}
Let $\,\la \P X\, ;\, X\in\C\ra\,$ be a hyperdoctrine and define
\[
\hom_X(P,Q) := \Pi_X(P\imp Q):(\P X)\op\tm\P X\to\P 1   
\]
\[
\meets_X(P,Q) := \Si_X(P\tm Q):\P X\tm\P X\to\P 1
\] 
where the quantification indexes denote the map $X\to 1$.
Then $\hom_X$ enriches $\P X$ in $\P 1$ and, for any map $f:X\to Y$,
the following adjunction-like laws hold: 
\[
\hom_X(f\ast Q, P) \iso \hom_Y(Q,\Pi_f P)  \qv  \hom_X(P,f\ast Q) \iso \hom_Y(\Si_f P,Q)
\]
\[
\meets_X(f\ast Q, P) \iso \meets_Y(Q,\Si_f P)  \qv  \meets_X(P,f\ast Q) \iso \meets_Y(\Si_f P,Q)
\]
\end{corollary}
\pf
Propositions~\ref{prod} and~\ref{left}, and the Frobenius law itself, give:
\[
\Pi_X(f\ast Q \imp P) \iso \Pi_Y\Pi_f(f\ast Q \imp P) \iso \Pi_Y(Q \imp \Pi_f P) 
\]
\[
\Pi_X(P\imp f\ast Q) \iso \Pi_Y\Pi_f(P \imp f\ast Q) \iso \Pi_Y(\Si_f P \imp Q) 
\]
\[
\Si_X(P\tm f\ast Q) \iso \Si_Y\Si_f(P \tm f\ast Q) \iso \Si_Y(\Si_f P \tm Q)
\]
\epf
Say that a map $f:X\to Y$ is ``surjective" if $\Si_f\top_X \iso \top_Y$, where $\top_X$ is a terminal 
object of $\P X$.
\begin{corollary}       \label{hyp2}
If $f:X\to Y$ is surjective map then
\[
\Pi_X(f\ast Q) \iso \Pi_Y Q   \qv   \Si_X(f\ast Q) \iso  \Si_Y Q
\]
\end{corollary}
\pf
For the first one we have:
\[
\Pi_X(f\ast Q) \iso \hom_X(\top_X,f\ast Q) \iso \hom_Y(\Si_f\top_X,Q) \iso \hom_Y(\top_Y,Q) \iso \Pi_Y Q 
\] 
The proof of the second one follows the same pattern:
\[
\Si_X(f\ast Q) \iso \meets_X(\top_X,f\ast Q) \iso \meets_Y(\Si_f\top_X,Q) \iso \meets_Y(\top_Y,Q) \iso \Si_Y Q 
\]
\epf
Note that for $\,\la \P X\, ;\, X\in\Set\ra\,$, Corollary~\ref{hyp2} becomes the fact that the inverse image
functor along a surjective mapping $f:X\to Y$ preserves non-emptyness and reflects maximality: 
if $P\subseteq Y$ is non-empty so it is $f^{-1}P$ and if $f^{-1}P = X$ then $P=Y$.

There are three canonical ways to get a ``truth value" in $\P 1$ from $P\in\P X$, namely
quantifications along $X:X\to 1$ and evaluation at a point $x:1\to X$: 
\[  \Pi_X P \qv \Si_X P  \qv  x\ast P  \]
In the above proposition, we have used the fact that quantifications along $X$ are ``represented" (by $\top_X$):
\[
\Pi_X P \iso \hom_X(\top_X,P) \qv \Si_X P \iso \meets_X(\top_X,P)
\]
Now we show that the same is true for evaluation; 
namely, evaluation at $x$ is ``represented" by the ``singleton":
\[
\{ x \} := \Si_x\top_1
\]   
\begin{corollary}       \label{hyp3}
Given a point $x:1\to X$ there are isomorphisms 
\[
x\ast P \iso \hom_X(\{x\},P)   \qv   x\ast P \iso \meets_X(\{x\},P)
\]
natural in $P\in\P X$.
\end{corollary}
\pf
\( \hom_X(\Si_x\top_1,P)  \iso \hom_1(\top_1,x\ast P)  \iso \Pi_1(x\ast P) \iso x\ast P \)  \\

\( \qq \meets_X(\Si_x\top_1,P)  \iso \meets_1(\top_1,x\ast P)  \iso \Si_1(x\ast P) \iso x\ast P \) \\

(Note that the last index $1$ is the identity on $1\in\C$.)
\epf
\begin{remark}   \label{rmkhyp}
Suppose that $\C$ has pullbacks, so that we also have the doctrine $\,\la \C/X\, ;\, X\in\C\ra\,$,
with $f_!\adj f\st:\C/Y\to\C/X$ for $f:X\to Y$.
Suppose also that $\,\la \P X\, ;\, X\in\C\ra\,$ satisfies the comprehension axiom~\cite{law70} $c_X \adj k_X:\P X\to\C/X$.
Then the set-valued ``external evaluation" of $P\in\P X$ at $x:1\to X$ can be expressed in various ways:
\[
\P 1(\top_1, x\st P) \iso \P X(\{x\},P) \iso \P X(c_X x,P) \iso
\]
\[
\iso \C/X(x,k_X P) \iso \C/X(x_!1,k_X P) \iso \C(1,x\st k_X P)
\]
\end{remark}
\begin{corollary}   [formulas for quantifications]    \label{hyp4}
Given $P\in\P X$, a map $f:X\to Y$ and a point $y:1\to Y$, there are isomorphisms 
\[
y\ast \Pi_f P \iso \hom_X(f\ast\{x\},P)   \qv   y\ast \Si_f P \iso \meets_X(f\ast\{x\},P)
\]
natural in $P\in\P X$.
\end{corollary}
\pf
\( y\ast \Pi_f P  \iso \hom_Y(\{y\},\Pi_f P)  \iso \hom_X(f\ast\{y\},P) \) \\

\(\qq y\ast \Si_f P  \iso \meets_Y(\{y\},\Si_f P)  \iso \meets_X(f\ast\{y\},P) \)
\epf
Note that for $\,\la \P X\, ;\, X\in\Set\ra\,$, Corollary~\ref{hyp4} gives the classical formula for the coimage
of a part along a mapping $f$, and a (less classical) formula for the image: 
$y$ is in the image $\Si_f P$ iff its inverse image meets $P$.

\section{Temporal doctrines}
\label{temp}

A {\bf temporal doctrine} $\,\la i_X:\M X\to\P X\ot\M'X:i'_X\, ;\, X\in\C\ra\,$ consists of two indexed functors
with the same codomain, satisfying the axioms listed below.

We denote the substitution functors along a map $f:X\to Y$ in $\C$ by 
\[ 
f\cd :\M Y\to\M X \qv f'\cd :\M'Y\to\M'X \qv f\ast :\P Y\to\P X
\]
Thus we have (coherent) isomorphisms:
\[
(gf)\cd \iso f\cd g\cd \qv (gf)'\cd \iso f'\cd g'\cd \qv (gf)\ast \iso f\ast g\ast
\]
(and similarly for identities) and also
\[
i_X f\cd \iso f\ast i_Y  \qv  i'_X f'\cd \iso f\ast i'_Y
\]
We denote by $\B X$ the indexed pullback $\M X \tm_{\P X} \M' X$, by $j_X$ and $j'_X$
its indexed projections to $\M X$ and $\M'X$ respectively, and
\[
b_X:= i_Xj_X = i'_Xj'_X : \B X \to \P X  
\] 
The first group of axioms requires the existence of some adjoint functors:
\begin{enumerate}
\item
The indexing category has a terminal object: $1\in\C$.
\item
The categories of $\P X$ are cartesian closed.
Thus, for any $X\in\C$, we have a terminal object $1_X\in\P X$, products $P\tm_X Q$ and
exponentials $P\imp_X Q$.
\item
The substitution functors $f\ast :\P Y\to\P X$ have both left and right adjoints:
\[
\Si_f \adj f\ast \adj \Pi_f 
\]  
\item
The functors $i_X:\M X\to\P X$ and $i'_X:\M'X\to\P X$ have both left and right adjoints:
\[
\dm_X \adj i_X \adj \,\sqq_X  \qv  \dm'_X \adj i'_X \adj \,\sqq'_X
\]  
\item
The doctrine $\P X$ satisfies the comprehension axiom~\cite{law70}:
the canonical functors $c_X:\C/X\to\P X$ (sending $f:T\to X$ to $\Si_f 1_T$) have right adjoints:
\[
c_X \adj k_X:\P X\to\C/X
\]
\end{enumerate}
The second group of axioms imposes some exactness condition on these functors:
\begin{enumerate}
\item
The functors $i_X$ and $i'_X$ are fully faithful:
\[
\dm_X i_X \iso \id_{\M X}  \qv \sq_X i_X \iso \id_{\M X}
\]
(and similarly for $i'_X$). 
\item
The doctrine $\P X$ satisfies the Frobenius law:
\eq        
\Si_f P \tm_Y Q \iso \Si_f (P\tm_X f\ast Q) 
\eeq
for any $f:X\to Y$ (naturally in $P\in\P X$ and $Q\in\P Y$).
\item
The adjunctions $\dm_X \adj i_X$ and $\dm'_X \adj i'_X$ satisfy the ``mixed Frobenius laws",
that is their units induce isomorphisms
\eq
\dm_X(P \tm_X i'_X N) \iso \dm_X(i_X\dm_X P \tm_X i'_X N)  \label{rsl1}
\eeq
\eq
\dm'_X(P \tm_X i_X M) \iso \dm'_X(i'_X\dm'_X P \tm_X i_X M)   \label{rsl2}
\eeq
(natural in $P\in\P X$, $N\in\M'X$ and $M\in\M X$).
\item
The projections $j_1:\B 1 \to \M 1$ and $j'_1:\B 1 \to \M'1$ are isomorphisms.
\item
The comprehension functors $k_X:\P X\to\C/X$ are fully faithful:
\eq
c_X k_X P = \Si_{k_X P}1_{X_!(k_X P)} \iso P
\eeq
(where we use the notations of Remark~\ref{rmkhyp}, so that $X_!$ is the domain projection $\C/X\to \C$.
Note that the index $k_X P$ of $\Si$ is an object of $C/X$, so that it should be more exactly 
be replaced by $X_!(k_X P)$, where now $k_X P$ denote the map to the terminal in $C/X$).
\end{enumerate}
\begin{examples}  \label{ex}
\begin{enumerate}
\item
Any hyperdoctrine $\,\la \P X\, ;\, X\in\C\ra\,$ (see Section~\ref{hyper}) with a fully faithful
comprehension functor gives rise to a (rather trivial) temporal doctrine:
\[ \la \,\id:\P X\to\P X\ot\P X:\id\, ;\, X\in\C\,\ra\, \]
Thus the results of Section~\ref{hyper} can be seen as particular cases of those we will obtain for temporal doctrines.
\item
$\,\la i_X:\D X\to\P X\ot\U X:i'_X\, ;\, X\in\Pos\ra\,$, where $\,\la \P X\, ;\, X\in\Pos\ra\,$ is the doctrine
of all the parts of a poset, while $\D X$ and $\U X$ are the subdoctrines of lower-closed and upper-closed parts of $X$.
\item
$\,\la i_X:\M X\to\Grph/X\ot\M'X:i'_X\, ;\, X\in\Grph\ra\,$, where $\Grph$ is the category of reflexive graphs,
while $\M X$ and $\M'X$ are the categories of left and right actions of $X$ (or of the free category generated by it).
\item
Groupoids or sets endowed with an equivalence relation give rise to ``symmetrical" temporal doctrines:
all the projections $j_X$ and $j'_X$ are isomorphisms.
Note that, since the axioms are symmetrical, each temporal doctrine has a dual obtained by exchanging the left
and the right side (that is $i$ and $i'$); while a symmetrical temporal doctrine is clearly self-dual
(that is isomorphic to its own dual) the same is true for $\Grph$, via the ``opposite" functor $\Grph \to \Grph$.
\item
Any strong balanced factorization category $\la \C\,;\, \E,\M \ra$ ~\cite{pis08,pis09} such that $\C$ is locally cartesian
closed and $\M/X$ and $\M'/X$ are coreflective in $\C/X$, gives rise to the temporal doctrine
\( \la i_X:\M/X\to\C/X\ot\M'/X:i'_X\, ;\, X\in\C\ra \)
\item
Given a temporal doctrine on a category $\C$ and any subcategory $\C'$ of $\C$ such that $1\in\C'$, one gets by restriction
another temporal doctrine on $\C'$. 
\end{enumerate}
\end{examples}

\begin{remark}
The name ``temporal doctrine" is clearly suggested by the functors $\dm$, $\sqq$, $\dm'$ and $\sqq'$,
which can be seen as modal operators acting in the two directions of time.

A categorical approach to modal and tense logic was developed in the eighties by 
Ghilardi and Meloni and indipendently by Reyes et al. 
Not being here specifically concerned with these logics, we just note that the temporal
doctrine of posets mentioned in the examples above is also an instance of temporal doctrine in the sense of~\cite{mel}.

Let me also acknoweldge that it was prof. Giancarlo Meloni, the supervisor of my phd thesis,
who introduced me to categorical logic showing in particular how adjunctions can be an effective
tool for doing calculations.
\end{remark}

\section{Basic properties}
\label{fp}

\subsection{terminology}
\label{term}

Since a (weaker form of) temporal doctrine is mainly intended to model the situation
$\la \Set^{X\op}\to\Cat/X\ot\Set^X\, ;\, X\in\Cat\ra$, the objects of $\C$ should be thought of
as generalized categories. In fact in the sequel we will freely borrow terminology from category
theory, whenever opportune. 
However, the interior $\sqq_X$ and closure $\dm_X$ operators suggest that 
it also make sense to consider the objects of $\C$ as a sort of spaces, 
so that we will also borrow some terminology from topology; in fact, the links with that subject can be taken 
quite seriously as sketched in~\cite{pis09}, where it is discussed also the significance of the ``closure" reflection
in ``open parts" (or ``local homeomorphisms").
Anyway, if $X$ is a topological space and $i_X$ and $i'_X$ are the inclusion of open and closed parts respectively in $\P X$,
the mixed Frobenius laws (and their equivalent ones)
hold true when they make sense, that is when only the operators $\sqq_X$ and $\dm'_X$ are involved.

Thus we sometimes refer to objects and arrows of $\C$ as ``spaces" and ``maps"; to the objects of $\P X$ 
as ``parts" of $X$ and to those of $\M X$ and $\M'X$ as left closed and right closed parts of $X$, respectively.
The reflections $\dm_X$ and $\dm'_X$ are the left and right ``closure" operators respectively,
while $\sqq_X$ and $\sqq'_X$ are the left and right ``interior" operators.

Apart from the axioms concerning the comprehension adjunctions, 
a temporal doctrine is a hyperdoctrine $\P X$ (in the sense of Section~\ref{hyper}) 
with two reflective and coreflective indexed subcategories $\M X$ and $\M'X$ 
such that $\M 1$ and $\M'1$ are isomorphic as subcategories of $\P 1$;
furthermore, and most importantly, we assume the mixed Frobenius laws~(\ref{rsl1}) and~(\ref{rsl2}),   
which are rich of important consequences. The reason of their name follows by Remark~\ref{ff}:
they look like the Frobenius laws for $\dm_X \adj i_X$ and  $\dm'_X \adj i'_X$, except that $i_X$ and $i'_X$
are exchanged in the second factors.
In fact, we have the corresponding mixed equivalent conditions:
\begin{prop}   \label{temp1}
The following laws hold in a temporal doctrine, and each o them can be used in the definition
in place of~{\rm (\ref{rsl1})} and~{\rm(\ref{rsl2})}:
\eq   \label{rsl3}
i'_X N\imp i_X M \iso i_X\sq_X(i'_X N\imp i_X M)  \,\,\,\, ; \,\,\,  i_X M\imp i'_X N \iso i'_X\sq'_X(i_X M\imp i'_X N)
\eeq
\eq    \label{rsl4}
\sq'_X(P \imp i_X M) \iso \sq'_X(i_X\dm_X P \imp i_X M)  \,\,\,\, ; \,\,\,   \sq_X(P \imp i'_X N) \iso \sq_X(i'_X\dm'_X P \imp i'_X N)
\eeq
Furthermore
\[
\sq_X(i_X M\imp P) \iso \sq_X(i_X M \imp i_X\sq_X P)  \,\,\,\, ; \,\,\,  \sq'_X(i'_X N\imp P) \iso \sq'_X(i'_X N \imp i'_X\sq'_X P)
\]
\end{prop} 
\pf
As in Proposition~\ref{left}, both the members of~(\ref{rsl1}) and of~(\ref{rsl2}) have two right adjoints,
one for each parameter considered, giving the conditions above.
For the last statement, recall~(\ref{ff1}). 
\epf
From~(\ref{rsl3}) we immediatley get:
\begin{corol}
If the part $P\in\P X$ is left closed and $Q\in\P X$ is right closed (that is $P\iso i_X M$ and $Q\iso i'_X N$)
then $P\imp Q$ is itself right closed.
\epf
\end{corol}
As already mentioned in the Introduction, we so have an ``explanation" of the fact that the complement 
of an upper-closed part of a poset is lower-closed (and conversely).
\begin{corol}
The categories $\B X$ are themselves cartesian closed, with the ``same" exponential of $\P X$: 
\[
b_X B \imp b_X C \iso b_X(B \imp_{\B X} C)
\] 
\epf
\end{corol}
\begin{prop}     \label{temp2}
$\,\la \M X\, ;\, X\in\C\ra\,$ and $\,\la \M'X\, ;\, X\in\C\ra\,$ are themselves hyperdoctrines, 
with a fully faithful comprehension adjoint.
\end{prop}
\pf
\begin{enumerate}
\item
As in Remark~\ref{ff}, the categories $\M X$ and $\M'X$ are cartesian closed,
with exponentials given by 
\[
\sq_X (i_X L \imp i_X M)    \qv   \sq'_X (i'_X N \imp i'_X O)
\] 
We denote products in $\M X$ and $\M'X$ by
\[
\top_X  \qv  L\wdg_X M  \qv  \top'_X  \qv  N\wdg'_X O
\]
\item
The substitution functors for left and right closed parts have both left and right adjoints:
\[
\ex_f \adj f\cd \adj \fa_f  \qv  \ex'_f \adj f'\cd \adj \fa\,'_f 
\] 
where 
\eq     \label{comm}
\ex_f \iso \dm_Y\Si_f\, i_X  \qv  \fa_f \iso \sq_Y\Pi_f\, i_X
\eeq
(and similarly for $\ex'_f$ and $\fa\,'_f$).
Note that these satisfy:
\eq   
\ex_f \dm_X \iso \dm_Y\Si_f   \qv  \fa_f \sq_X \iso \sq_Y\Pi_f
\eeq  
(and similarly for $\ex'_f$ and $\fa\,'_f$).
\item
The canonical functors $\C/X\to\M X$ send $f:T\to X$ to 
\[
\ex_f \top_T \iso \dm_X\Si_f\, i_X\top_T \iso \dm_X\Si_f 1_T \iso \dm_X c_X f
\]
that is factor through the corresponding ones for $\P X$.
Thus, they have the functors $k_X i_X:\M X\to\C/X$ as fully faithful right adjoints
(and similarly for $\M'X$; we leave it to the reader to check the above factorization for the arrow mapping).
\end{enumerate}
The fact that the adjunctions $\dm_X\adj i_X$ and $\dm'_X\adj i'_X$ satisfy the mixed Frobenius laws
implies a restricted form of the Frobenius law for each of them and also
for $\ex_f\adj f\cd$ and $\ex'_f\adj f'\cd$, which will be used in the sequel:
\begin{prop}    [restricted Frobenius laws] \label{temp3}
For any $X\in\C$, there are natural isomorphisms:
\[
\dm_X(P \tm_X i_X j_X B) \iso \dm_X P \wdg_X j_X B
\]
For any $f:X\to Y$ in $\C$, there are natural isomorphisms:
\[
\ex_f(M \wdg_X f\cd j_Y B) \iso \ex_f M \wdg_X j_Y B  \qv  \ex'_f(N \wdg'_X f'\cd j'_Y B) \iso \ex'_f N \wdg'_X j'_Y B 
\]
\end{prop}
\pf
For the first one, by the mixed Frobenius law we get:
\[
\dm_X(P \tm_X i_X j_X B) \iso \dm_X(P \tm_X i'_X j'_X B) \iso \dm_X(i_X\dm_X P \tm_X i'_X j'_X B) \iso 
\]
\[
\iso \dm_X(i_X\dm_X P \tm_X i_X j_X B) \iso \dm_X i_X(\dm_X P \wdg_X j_X B) \iso \dm_X P \wdg_X j_X B
\]
For the second one, we then have:
\[   
\ex_f(M \wdg_X f\cd j_Y B)  \iso \dm_Y\Si_f i_X (M \wdg_X f\cd j_Y B) \iso 
\]
\[
\iso \dm_Y\Si_f(i_X M \tm_X  i_X f\cd j_Y B) \iso \dm_Y\Si_f(i_X M \tm_X f\ast i_Y j_Y B)  \iso
\]
\[
\iso \dm_Y(\Si_f i_X M \tm_Y i_Yj_Y B) \iso (\dm_Y\Si_f i_X M) \wdg_Y j_Y B  \iso \ex_f M \wdg_Y j_Y B 
\]
\epf

\section{Functors valued in truth values}
\label{tv}

In the sequel, a major role will be played by the ``truth values" category $\B 1$.
We denote by $\true$ its terminal object, so that 
\[
j_1\true \iso \top_1  \qv  j'_1\true \iso \top'_1
\]
The functors $X\ast b_1:\B 1\to\P X$ can be factorized in various ways:
\[
i_X X\cd j_1 = X\ast i_1j_1 = X\ast i'_1j'_1 = i'_X X'\cd j'_1
\]
(where $X$ denotes also the map $X\to 1$).
Thus their left and right adjoints can be factorized as:
\eq    \label{coend}
j_1^{-1}\ex_X\dm_X \iso j_1^{-1}\dm_1\Si_X \iso {j'_1}^{-1}\dm'_1\Si_X \iso {j'_1}^{-1}\ex\,'_X\dm'_X 
\eeq
\eq    \label{end}
j_1^{-1}\fa_X\sq_X \iso j_1^{-1}\sq_1\Pi_X \iso {j'_1}^{-1}\sq'_1\Pi_X \iso {j'_1}^{-1}\fa\,'_X\sq'_X 
\eeq
We refer to (anyone of) these as the ``coend" and ``end" functors at $X$, respectively:
\[
\coend_X\adj X\ast b_1 \adj\eend_X : \P X \to \B 1
\]
This terminology is justified by the fact that, for a bifunctor $H : X\op\tm X \to \Set$, one can
easily construct an object $h$ of $\Cat/X$ such that $\eend_X h$ gives the usual end of $H$, while
$\coend_X h$ gives the coend of $H$ in the sense of strong dinaturality, which in most relevant
cases reduces to the usual one as well (see~\cite{pis07}).

Next we define the functors 
\[
\meets_X : \P X \tm \P X \to \B 1  \qv  \hom_X : (\P X)\op \tm \P X \to \B 1
\]
\[
\meets_X(P,Q) := \coend_X(P\tm Q)  \qv  \hom_X(P,Q) := \eend_X(P\imp Q)
\]
and their restrictions
\[
\ten_X : \M'X \tm \M X \to \B 1   
\]
\[
\nat_X : (\M X)\op \tm \M X \to \B 1  \qv  \nat'_X : (\M'X)\op \tm \M'X \to \B 1
\]
\[
\ten_X(N,M) := \meets_X(i'_X N,i_X M)    
\]
\[
\nat_X(L,M) := \hom_X(i_X L,i_X M)  \qv  \nat'_X(N,O) := \hom_X(i'_X N,i'_X O)
\]
For instance, in the temporal doctrine of posets $\B 1 \iso \P 1 \iso \{\true,\false\}$ and 
\[ \meets_X(P,Q) = \true \] 
iff $P$ and $Q$ have a non-empty intersection (and similarly for $\ten_X(N,M)$).
Of course, $\hom_X(P,Q) = \true$ iff $P\subseteq Q$ (and similarly for $\nat_X(L,M)$ and $\nat'_X(N,O)$).
In the temporal doctrine of reflexive graphs, $\P 1 \iso \Grph$ while $\B 1 \iso \Set$.

Note that $\hom_X$ and $\ten_X$ are valued in $\B 1$ rather than in $\P 1$ as in Section~\ref{hyper},
so that the notation is in fact consistent only for the first example in~\ref{ex}.

\subsection{The enriched ``adjunction" laws}

In the following proposition, we show that the adjunctions which define a temporal doctrine 
can be internalized, that is they are enriched in the truth values category $\B 1$.
Furthermore, some of them have an exact counterpart in a similar law, with
the ``meets" or ``tensor" functors in place of the ``hom" or ``nat" functors;
the proofs are also nicely symmetrical.  
\begin{prop}      \label{td}
The functors $\hom_X$, $\nat_X$ and $\nat'_X$ enrich $\P X$, $\M X$ and $\M'X$
respectively in $\B 1$ and, for any space $X\in\C$ or map $f:X\to Y$,
there are natural isomorphisms: 
\eq   \label{adj1}
\hom_X(f\ast Q, P) \iso \hom_Y(Q,\Pi_f P)  \qv  \hom_X(P,f\ast Q) \iso \hom_Y(\Si_f P,Q)
\eeq
\eq   \label{adj2}
\meets_X(f\ast Q, P) \iso \meets_Y(Q,\Si_f P) \,\,\,\,\,;\,\,\,\,\,  \meets_X(P,f\ast Q) \iso \meets_Y(\Si_f P,Q)
\eeq
\eq   \label{adj3}
\nat_X(M,\sq_X P) \iso \hom_X(i_X M,P)  \qv  \nat'_X(N,\sq'_X P) \iso \hom_X(i'_X N,P) 
\eeq
\eq   \label{adj4}
\nat_X(\dm_X P, M) \iso \hom_X(P,i_X M)  \qv  \nat'_X(\dm'_X P, N) \iso \hom_X(P,i'_X N) 
\eeq
\eq   \label{adj5}
\ten_X(N,\dm_X P) \iso \meets_X(i'_X N,P)  \qv  \ten_X(\dm'_X P,M) \iso \meets_X(P,i_X M) 
\eeq
\eq   \label{adj6}
\nat_X(f\cd M, L) \iso \nat_Y(M,\fa_f L)  \qv  \nat'_X(f'\cd O, N) \iso \nat'_Y(O,\fa'_f N)
\eeq
\eq   \label{adj7}
\nat_X(L,f\cd M) \iso \nat_Y(\ex_f L,M)  \qv  \nat'_X(N,f'\cd O) \iso \nat'_Y(\ex'_f N,O)
\eeq
\eq   \label{adj8}
\ten_X(f'\cd N,L) \iso \ten_Y(N,\ex_f L)  \qv  \ten_X(N,f\cd M) \iso \ten_Y(\ex'_f N,M)
\eeq
\end{prop}
\pf
For the first part, see Proposition~\ref{prod}.
For~(\ref{adj1}) and~(\ref{adj2}), recall Corollary~\ref{hyp} and note that the present 
``hom" and ``meets" functors factor through the ones there.

Equations~(\ref{adj3}),~(\ref{adj4}) and~(\ref{adj5}) follow from Proposition~\ref{temp1} and
the other factorizations of the coend and the end functors in~(\ref{coend}) and~(\ref{end}).
Recalling~(\ref{comm}), we obtain the remaining ones by composition of (enriched) adjoints.
Alternatively, one can explicitly derive them as we exemplify for~(\ref{adj7}):
\[
j_1^{-1}\sq_1\Pi_X(i_X L \imp i_X f\cd M) \iso j_1^{-1}\sq_1\Pi_Y\Pi_f(i_X L \imp f\ast i_Y M) 
\iso j_1^{-1}\sq_1\Pi_Y(\Si_f i_X L \imp i_Y M) 
\] 
\[
\iso {j'}_1^{-1}\fa\,'_Y\sq'_Y(\Si_f i_X L \imp i_Y M) \iso {j'}_1^{-1}\fa\,'_Y\sq'_Y(i_Y \dm_Y \Si_f i_X L \imp i_Y M)
\]
\epf

\section{Limits, colimits and Yoneda properties}
\label{lim}

As we will see in Section~\ref{cat}, most of the laws in Proposition~\ref{td} (namely those not containing $\hom$)
still hold for weak temporal doctrines, which include the motivating instance 
$\la \Set^{X\op}\to\Cat/X\ot\Set^X\, ;\, X\in\Cat\ra$.
Thus, with the same technique exploited in Section~\ref{hyper}, we begin to draw some consequences 
which in fact hold in the weaker context as well.
Accordingly, we mainly maintain the policy of using terms which reflect the case $\C = \Cat$ just mentioned.  

We define the (internal) ``limit" and ``colimit" functors by restricting the end and the coend functors
to (left or right) closed parts:
\[ \lim_X := \eend_X i_X \iso j_1^{-1}\fa_X : \M X \to \B 1  \] 
\[ \lim'_X := \eend_X i'_X \iso {j'_1}^{-1}\fa'_X : \M'X \to \B 1 \]
\[ \colim_X := \coend_X i_X \iso j_1^{-1}\ex_X : \M X \to \B 1  \]  
\[ \colim'_X := \coend_X i'_X \iso {j'_1}^{-1}\ex'_X : \M'X \to \B 1 \]

\begin{corollary}     [final maps preserve limits]  \label{td2}
Let $f:X\to Y$ be a final map, that is $\ex_f\top_X \iso \top_Y$. Then
\eq    \label{plim}
\lim_X(f\cd M) \iso \lim_Y M   \qv   \colim'_X(f'\cd N) \iso  \colim'_Y N
\eeq
Dually, if $f:X\to Y$ is initial, that is $\ex'_f\top'_X \iso \top'_Y$, then
\[
\colim_X(f\cd M) \iso \colim_Y M   \qv   \lim'_X(f'\cd N) \iso  \lim'_Y N
\]
\end{corollary}
\pf
For the first one of~(\ref{plim}), using~(\ref{adj7}) we have:
\[
\lim_X(f\cd M) \iso \nat_X(\top_X,f\cd M) \iso \nat_Y(\ex_f\top_X,M) \iso \nat_Y(\top_Y,M) \iso \lim_Y M 
\] 
For the second one of~(\ref{plim}), we follow exactly the same pattern using~(\ref{adj8}) instead:
\[
\colim'_X(f'\cd N) \iso \ten_X(f'\cd N,\top_X) \iso \ten_Y(N,\ex_f\top_X) \iso \ten_Y(N,\top_Y) \iso \colim'_Y N
\]
The other ones are proved in the same way.
\epf

There are three canonical ways of obtaining a truth value in $\B 1$ from a closed part, namely
the limit or the colimit functors and ``evaluation" at a point $x:1\to X$: 
\[  \lim_X M \qv \colim_X M  \qv  j_1^{-1}(x\cd M)  \]
\[  \lim'_X N \qv \colim'_X N  \qv  {j'_1}^{-1}(x'\cd N)  \]
In the above proposition, we have used the fact that limits and colimits over $X$ are ``represented" (by $\top_X$ or $\top'_X$).
Now we show that the same is true for evaluation; 
namely, evaluation at $x$ is ``represented" by the left and right ``slices":
\eq   \label{tde}
X/x := \ex_x \top_1   \qv   x\bs X := \ex'_x \top'_1 
\eeq   
Note that slices can be obtained as the (left or right) closure of the ``singletons" 
$\{x\} = \Si_x 1_1 = c_X x$ (see Section~\ref{hyper}):
\[
X/x = \dm_X \{x\}  \qv  x\bs X = \dm'_X \{x\}
\]
\begin{corollary}     [Yoneda properties]    \label{td3}
Given a point $x:1\to X$ in $\C$, there are isomorphisms 
\[
j_1^{-1}(x\cd M) \iso \nat_X(X/x,M)   \qv   j_1^{-1}(x\cd M) \iso \ten_X(x\bs X,M)
\]
natural in $P\in\P X$ (and dually for right closed parts).
\end{corollary}
\pf
\[
\nat_X(\ex_x\top_1,M) \iso \nat_1(\top_1,x\cd M) \iso \lim_1(x\cd M) \iso j_1^{-1}\fa_1(x\cd M) \iso j_1^{-1}(x\cd M)
\]
\[
\ten_X(\ex'_x\top'_1,M)  \iso \ten_1(\top'_1,x\cd M) \iso \colim_1(x\cd M) \iso j_1^{-1}\ex_1(x\cd M) \iso j_1^{-1}(x\cd M)
\]
\epf
\begin{corollary}    [formulas for quantifications, interior and closure]   \label{td4}
Given a map $f:X\to Y$ and a point $y:1\to Y$ in $\C$, there are isomorphisms 
\[
j_1^{-1}(y\cd\fa_f M) \iso \nat_X(f\cd Y/y,M)   \qv   j_1^{-1}(y\cd\ex_f M) \iso \ten_X(f'\cd y\bs Y,M)
\]
natural in $M\in\M X$ (and dually for right closed parts). 
There are isomorphisms 
\[
j_1^{-1}(x\cd\sq_X P) \iso \hom_X(i_X X/x,P)  \qv  j_1^{-1}(x\cd\dm_X P) \iso \meets_X(i'_X x\bs X,P)
\]
natural in $P\in\P X$ (and dually for right closure). 
\end{corollary}
\pf
Using Corollary~\ref{td3} and~(\ref{adj6}),~(\ref{adj8}),~(\ref{adj3}) and~(\ref{adj5}) respectively, we get:

\(
j_1^{-1}(y\cd \fa_f M)  \iso \nat_Y(Y/y,\fa_f M) \iso \nat_X(f\cd Y/y,M)
\)

\(
j_1^{-1}(y\cd \ex_f M)  \iso \ten_Y(y\bs Y,\ex_f M) \iso \ten_X(f'\cd y\bs Y,M)
\)

\(
j_1^{-1}(x\cd\sq_X P) \iso \nat_X(X/x,\sq_X P) \iso \hom_X(i_X X/x,P)  
\)

\(
j_1^{-1}(x\cd\dm_X P) \iso \ten_X(x\bs X,\dm_X P) \iso \meets_X(i'_X x\bs X,P)
\)
\epf

\section{Exploiting comprehension}
\label{compr}

In this section and in the next one we present some consequences of the comprehension 
adjunction $c_X\adj k_X:\P X\to\C/X$ and of the assumption that it is fully faithful.

\subsection{The components functor}

We define the ``components" functor $\comp:\C\to\B 1$ by:
\[
\comp X := \coend_X 1_X = \colim_X \top_X = \colim'_X \top'_X
\]
\begin{remarks}  \label{con}
Note that $X$ is connected, that is $\comp X \iso \true$, iff $X\to 1$ is final (or initial).
Note also that the components functor
\[
\comp X = \coend_X 1_X = j_1^{-1}\dm_1\Si_X 1_X = j_1^{-1}\dm_1 c_1 X
\]
is left adjoint to the full inclusion $k_1 i_1 j_1 : \B 1 \to \C$.
Coherently, we say that a space $X\in\C$ is ``discrete" if $X\iso k_1 i_1 j_1 V$, for a truth value $V\in\B 1$,
so that $\comp$ yields in fact the reflection in discrete spaces.
\end{remarks}
Conversely, the coend functor can be reduced to components or to a colimit by
\eq      \label{compr0}
\coend_X P \iso \comp X_! k_X P \iso \colim_{X_! k_X P}\top_{X_! k_X P} 
\eeq
Indeed we have:
\[
\coend_X P \iso j_1^{-1}\dm_1\Si_X P \iso j_1^{-1}\dm_1\Si_X c_X k_X P \iso j_1^{-1}\dm_1\Si_X\Si_{k_X P} 1_{X_! k_X P} \iso
\]
\[
\iso j_1^{-1}\dm_1\Si_{X_! k_X P} 1_{X_! k_X P} \iso \coend_{X_! k_X P} 1_{X_! k_X P} \iso \comp X_! k_X P
\]

\subsection{The limit and colimit formulas for nat and ten}

Since $\coend_X$ can be reduced to a colimit by~(\ref{compr0}), the same is true for the functors $\meets_X$ and $\ten_X$.
In fact we can do better, by reducing $\ten_X(N,M)$ to a colimit
over the space corresponding to $i'_X N$ (or to $i_X M$), rather than to $i'_X N\tm_X i_X M$;
similarly, $\nat_X$ can be reduced to a limit:
\begin{prop}  [(co)limit formulas for ten and nat]  \label{tennat}
\eq    \label{compr1}
\meets_X(P,Q) \iso \coend_{X_! k_X Q}(k_X Q)\ast P \iso \coend_{X_! k_X P}(k_X P)\ast Q
\eeq
\eq    \label{compr2}
\ten_X(N,M) \iso \colim_{X_! k_X i'_X N}(k_X i'_X N)\cd M \iso \colim'_{X_! k_X i_X M}(k_X i_X M)'\cd N
\eeq
\eq     \label{compr3}
\hom_X(P,Q) \iso \eend_{X_! k_X Q}(k_X Q)\ast P \iso \eend_{X_! k_X P}(k_X P)\ast Q
\eeq
\eq     \label{compr4}
\nat_X(L,M) \iso \lim_{X_! k_X i_X L}(k_X i_X L)\cd M  \,\,\,\, ; \,\,\,\,  \nat'_X(N,O) \iso \lim'_{X_! k_X i'_X N}(k_X i'_X N)\cd O
\eeq
\end{prop}
\pf
For~(\ref{compr1}), by applying~(\ref{adj2}) of Proposition~\ref{td} we get
\[
\meets_X(P,Q) \iso \meets_X(P,\Si_{k_X Q}1_{X_!k_X Q}) \iso 
\]
\[
\iso \meets_{X_!k_X Q}((k_X Q)\ast P,1_{X_!k_X Q}) \iso \coend_{X_! k_X Q}(k_X Q)\ast P
\]
For~(\ref{compr2}), we then have:
\[
\ten_X(N,M) \iso \meets_X(i'_X N,i_X M) \iso \coend_{X_! k_X i'_X N}(k_X i'_X N)\ast i_X M \iso
\]
\[
\iso \coend_{X_! k_X i'_X N} i_{X_! k_X i'_X N} (k_X i'_X N)\cd M \iso \colim_{X_! k_X i'_X N}(k_X i'_X N)\cd M
\]
Similarly for~(\ref{compr3}), by applying the second of~(\ref{adj1}) instead we get:
\[
\hom_X(P,Q) \iso \hom_X(\Si_{k_X P}1_{X_!k_X P},Q) \iso 
\]
\[
\iso \hom_{X_!k_X P}(1_{X_!k_X P}, (k_X P)\ast Q) \iso \eend_{X_! k_X P}(k_X P)\ast Q
\]
Finally, for~(\ref{compr4}) we have:
\[
\nat_X(L,M) \iso \hom_X(i_X L,i_X M) \iso \eend_{X_! k_X i_X L}(k_X i_X L)\ast i_X M \iso
\]
\[
\iso \eend_{X_! k_X i_X L} i_{X_! k_X i_X L} (k_X i_X L)\cd M \iso \lim_{X_! k_X i_X L}(k_X i_X L)\cd M
\]
\epf
By applying~(\ref{compr2}) or~(\ref{compr4}) to the formulas for quantifications along 
a map $f:X\to Y$ of Corollary~\ref{td4}, we obtain a colimit and a limit formula for 
evaluation at $y:1\to Y$ of $\ex_f M$ (or $\ex'_f N$) and $\fa_f M$ (or $\fa'_f N$), respectively.

Say that a part $P\in\P X$ is ``left dense" if its left closure is terminal: $\dm_X P \iso \top_X$;
right density is of course defined dually. 
\begin{prop}    \label{compr5}
A part $P\in\P X$ is left dense iff $k_X P$ is final.
A map $f:X\to Y$ in $\C$ is final iff $c_Y f$ is left dense in $Y$.
A space $X\in\C$ is connected iff $\Si_X 1_X$ is dense.
\end{prop}
\pf
(Note that we implicitly use the canonical bijection between the objects of $\C/X$ and maps in $\C$ with
codomain $X$.)
For the first one we have:
\[
\dm_X P \iso \dm_X \Si_{k_X P}1_{X_!k_X P} \iso \dm_X \Si_{k_X P}i_{X_!k_X P}\top_{X_!k_X P} \iso \ex_{k_X P}\top_{X_!k_X P}
\]
and for the second one:
\( \qq \ex_f\top_X \iso \dm_Y \Si_f 1_X \iso \dm_Y c_Y f \)

For the last one, see Remarks~\ref{con}.
\epf

\section{The sup and inf reflections}
\label{sup}

For $X\in\C$, let $\ov X$ (resp. ${\ov X}'$) be the full subcategory of $\M X$ (resp. $\M'X$)
generated by the left (resp. right) slices~(\ref{tde}), and denote the inclusion functors by 
\[
h_X : \ov X \to \M X  \qv  h'_X : {\ov X}' \to \M'X
\]
The partially defined left adjoints to $i_X h_X$, $i'_X h'_X$, $k_X i_X h_X$ and $k_X i'_X h'_X$
are denoted respectively by:
\[
\sup_X:\P X \to \ov X \,\,\,\, ;\,\,\,\,  \inf_X:\P X \to {\ov X}'  \qv  
\Sup_X:\C/X \to \ov X  \,\,\,\, ;\,\,\,\,  \Inf_X:\C/X \to {\ov X}'
\]
\begin{remark}  \label{rmksup}
Of course, $\Sup_X f$ exists iff $\sup_X c_X f$ does, and in that case they are the same.
Note also that the sup (resp. inf) of a part depends only on its left (resp. right) closure: $\sup_X P$ exists 
iff $\dm_X P$ has a reflection in $\ov X$. 
\end{remark}
\begin{prop}
The following are equivalent for a space $X\in\C$:
\begin{enumerate}
\item
$X$ has a final point $x:1\to X$;
\item
there is a left dense part of $X$ with a sup;
\item
there is a final map $T\to X$ with a Sup;
\item
$\id_X:X\to X$ has a Sup.
\end{enumerate}
\end{prop}
\pf
By Proposition~\ref{compr5} and Remark~\ref{rmksup}, the last three conditions are equivalent.  
Since any point $x:1\to X$ has the slice $X/x$ as its Sup, 1) implies 3). 

Suppose conversely that a left dense part $P\in\P X$ has a sup; then by Remark~\ref{rmksup} 
$\dm_X P \iso \top_X$ has a reflection $X/x$ in $\ov X$.
By general well-known facts about reflections, it follows that also $X/x \iso \top_X$,
that is the point $x:1\to X$ is final.
\epf
\begin{prop}  [final maps preserve sups]
Let $t:S\to T$ be a final map in $\C$;
then $f:T\to X$ has a $\Sup$ iff $ft:S\to X$ does, and in that case they coincide.
\end{prop}
\pf
By Remark~\ref{rmksup}, it is enough to show that $\dm_X c_X ft \iso \dm_X c_X f$:
\[
\dm_X c_X ft \iso \ex_{ft}\top_S \iso \ex_f\ex_t\top_S \iso \ex_f\top_T \iso \dm_X c_X f
\]
\epf
Note that the sup (resp. inf) reflections give, for the weak temporal doctrine of categories, 
the colimit (resp. limit) of the corresponding functor
(see also~\cite{pis07} and~\cite{pis08}).  
Thus the above proposition can be seen as the ``external" correspective of Corollary~\ref{td2}.

\section{The logic of categories}
\label{cat}

Now we weaken the axioms of temporal doctrine so to include the motivating instance 
\[\la \Set^{X\op}\to\Cat/X\ot\Set^X\, ;\, X\in\Cat\ra \]
and show that most of the laws in Proposition~\ref{td} still hold.

\subsection{Weak temporal doctrines}

Weak temporal doctrines are defined as temporal doctrines except that:
\begin{enumerate}
\item
We do not require the existence of the $\Pi_f$, right adjoint to $f\ast$;
instead we do require, for any map $f:X\to Y$, the existence of 
\[
\fa_f:\M X\to \M Y \qv  \fa'_f:\M'X\to \M'Y
\]
right adjoint to $f\cd$ and $f'\cd$ respectively.
\item
We do not assume that the categories of parts $\P X$ are cartesian closed, but only 
that they are cartesian and that the left or right closed parts are exponentiable therein, that is
\[
i_X M\imp P  \qv  i'_X N\imp P
\]
always exist in $\P X$.
\end{enumerate}
In a weak temporal doctrine we still have the $\coend_X$ and $\eend_X$ functors as in Section~\ref{tv},
left and right adjoint to:
\[
i_X X\cd j_1 = X\ast i_1j_1 = X\ast i'_1j'_1 = i'_X X'\cd j'_1
\]
with the difference that, among the factorizatons of $\eend_X$ in~(\ref{end}), only 
\eq    \label{wend}
j_1^{-1}\fa_X\sq_X \iso {j'_1}^{-1}\fa\,'_X\sq'_X : \P X \to \B 1
\eeq
are always available.
We also define 
\[
\meets_X : \P X \tm \P X \to \B 1  \qv  \hom_X : (\P X)\op \tm \P X \to \B 1
\]
\[
\ten_X : \M'X \tm \M X \to \B 1  
\]
\[
\nat_X : (\M X)\op \tm \M X \to \B 1  \qv  \nat'_X : (\M'X)\op \tm \M'X \to \B 1
\]
as in Section~\ref{tv}, but now the $\hom_X$ may be only partially defined,
depending on whether the exponential $P\imp Q$ exists or does not exist.  
On the contrary, the above axiom on exponentials assures that the $\nat_X$ are always defined.
\begin{prop}      \label{wtd}
The functors $\nat_X$ and $\nat'_X$ enrich $\M X$ and $\M'X$
respectively in $\B 1$ and, for any space $X\in\C$ or map $f:X\to Y$,
there are natural isomorphisms: 
\eq   \label{wadj2}
\meets_X(f\ast Q, P) \iso \meets_Y(Q,\Si_f P) \,\,\,\,\,;\,\,\,\,\,  \meets_X(P,f\ast Q) \iso \meets_Y(\Si_f P,Q)
\eeq
\eq   \label{wadj5}
\ten_X(N,\dm_X P) \iso \meets_X(i'_X N,P)  \qv  \ten_X(\dm'_X P,M) \iso \meets_X(P,i_X M) 
\eeq
\eq   \label{wadj6}
\nat_X(f\cd M, L) \iso \nat_Y(M,\fa_f L)  \qv  \nat'_X(f'\cd O, N) \iso \nat'_Y(O,\fa'_f N)
\eeq
\eq   \label{wadj7}
\nat_X(L,f\cd M) \iso \nat_Y(\ex_f L,M)  \qv  \nat'_X(N,f'\cd O) \iso \nat'_Y(\ex'_f N,O)
\eeq
\eq   \label{wadj8}
\ten_X(f'\cd N,L) \iso \ten_Y(N,\ex_f L)  \qv  \ten_X(N,f\cd M) \iso \ten_Y(\ex'_f N,M)
\eeq
\end{prop}
\pf
Of course, the proofs of~(\ref{wadj2}),~(\ref{wadj5}) and~(\ref{wadj8}) are as in Proposition~\ref{td}.
For~(\ref{wadj6}), by applying Proposition~\ref{prod} to $f\cd\adj\fa_f$ and recalling that
exponentials in $\M X$ are given by $\sq_X(i_X L\imp i_X M)$, we have
\[
\nat_X(f\cd M, L) \iso \eend_X(i_X f\cd M \imp i_X L) \iso j_1^{-1}\fa_X\sq_X(i_X f\cd M \imp i_X L) \iso 
\]
\[
\iso j_1^{-1}\fa_Y\fa_f\sq_X(i_X f\cd M \imp i_X L) \iso j_1^{-1}\fa_Y\sq_Y(i_Y M \imp i_Y \fa_f L) \iso \nat_Y(M,\fa_f L)
\]  
For~(\ref{wadj7}), since we have adjunctions (with parameter $L\in\M X$)
\[
\ex_f(L\wdg_X X\cd j_1 V) \adj \nat_X(L,f\cd M)    \qv    \ex_f L \wdg_Y Y\cd j_1 V \adj \nat_Y(\ex_f L,M)
\]
and since, by the restricted Frobenius law of Proposition~\ref{temp3},
\[
\ex_f(L\wdg_X X\cd j_1 V) \iso \ex_f(L\wdg_X f\cd Y\cd j_1 V) \iso \ex_f L \wdg_Y Y\cd j_1 V 
\] 
the result follows from Lemma~\ref{par}.
\epf
\begin{prop}      \label{wtd2}
Whenever they are defined, there are natural isomorphisms: 
\eq   \label{wadj1}
\hom_X(f\ast Q, P) \iso \hom_Y(Q,\Pi_f P)  \qv  \hom_X(P,f\ast Q) \iso \hom_Y(\Si_f P,Q)
\eeq
\eq   \label{wadj3}
\nat_X(M,\sq_X P) \iso \hom_X(i_X M,P)  \qv  \nat'_X(N,\sq'_X P) \iso \hom_X(i'_X N,P) 
\eeq
\eq   \label{wadj4}
\nat_X(\dm_X P, M) \iso \hom_X(P,i_X M)  \qv  \nat'_X(\dm'_X P, N) \iso \hom_X(P,i'_X N) 
\eeq
\end{prop}
\pf
All of them can be proved as~(\ref{wadj7}) by showing that their left adjoints are isomorphic;
this requires the Frobenius law for the second of~(\ref{wadj1}) and 
the restricted Frobenius law for~(\ref{wadj4}).
We leave them as an exercise to the reader.
\epf
Since all we have proved in section~\ref{lim}, \ref{compr} and~\ref{sup} (except for~(\ref{compr3}) 
of Proposition~\ref{tennat}, wherein $\hom_X$ is only partially defined) depends only on the laws 
in Proposition~\ref{wtd} and on the comprehension axiom, those results hold in any weak temporal doctrine.
In particular we get, for ``generalized categories", the Yoneda properties, the formulas for quantifications 
(or ``Kan extensions"), the formulas for the (co)reflection in (left or right) ``closed parts" (or ``actions"),  
and the properties of final or initial maps with respect to (co)limits (both ``externally" and ``internally").

Thus we mantain that the logic of weak temporal doctrines well deserves to be called    
``a logic for categories", in fact
\begin{enumerate}
\item
being summarized by a few adjunction-like laws, it lends itself to effective and transparent calculations;
furthermore, along with the obvious ``left-right" symmetry, there is a far more interesting sort of duality:
 in many cases laws and proofs on the ``hom-side" correspond exactly to those on the ``tensor-side"; 
\item
this calculus allows one to easily derive some basic non-trivial categorical facts;
\item
it is ``autonomous", providing its own truth values;
\item
suitable natural strenghtenings or weakenings can be considered, so to obtain more refined properties
or a wider range of applications; some of them will be considered in a forthcoming work.
\end{enumerate}

\begin{refs}

\bibitem[Ghilardi \& Meloni, 1991]{mel} S. Ghilardi, G. Meloni (1991), {\em Relational and Topological Semantics
for Temporal and Modal Predicative Logics}, Atti del Congresso: Nuovi problemi della Logica e della Filosofia della Scienza
(Viareggio), CLUEB Bologna, 59-77.

\bibitem[Johnstone, 2002]{ele} P.T. Johnstone (2002), Sketches of an Elephant: A Topos Theory Compendium,
Oxford Science Publications.

\bibitem[Kan, 1958]{kan} D.M. Kan (1958), {\em Adjoint Functors},
{Trans. Amer. Math. Soc.} {\bf 87}, 294-329. 
 
\bibitem[Lawvere, 1970]{law70} F.W. Lawvere (1970), {\em Equality in Hyperdoctrines and the Comprehension Scheme as an Adjoint Functor},
Proceedings of the AMS Symposium on Pure Mathematics, XVII, 1-14. 

\bibitem[Lawvere, 1989]{law89} F.W. Lawvere (1989), {\em Qualitative Distinctions between some Toposes of Generalized Graphs},
Proceedings of the AMS Symposium on Categories in Computer Science and Logic, Contemporary Mathematics, vol. 92, 261-299.

\bibitem[Pisani, 2007]{pis07} C. Pisani (2007), Components, Complements and the Reflection Formula, 
{\em Theory and Appl. Cat.} {\bf 19}, 19-40. 

\bibitem[Pisani, 2008]{pis08} C. Pisani (2008), Balanced Category Theory, 
{\em Theory and Appl. Cat.} {\bf 20}, 85-115. 

\bibitem[Pisani, 2009]{pis09} C. Pisani (2009), Balanced Category Theory II, preprint, arXiv:math.CT/0904.1790v3.

\end{refs}

\end{document}